\documentclass[11pt,leqno]{article}
\usepackage{amsmath,amsfonts,amscd,amssymb,theorem}

\long\def\comment#1\endcomment{}

\makeatletter
\begingroup
\gdef\th@dotted{\normalfont\itshape
  \def\@begintheorem##1##2{%
        \item[\hskip\labelsep \theorem@headerfont ##1\ ##2.]}%
\def\@opargbegintheorem##1##2##3{%
   \item[\hskip\labelsep \theorem@headerfont ##1\ ##2\ (##3).]}}
\endgroup
\makeatother

\theoremstyle{dotted}

\newtheorem{theorem}{Theorem}[section]
\newtheorem{lemma}[theorem]{Lemma}

\newtheorem{prop}[theorem]{Proposition}
\newtheorem{corr}[theorem]{Corollary}

\makeatletter
\begingroup
\gdef\th@upshape{\normalfont
  \def\@begintheorem##1##2{%
        \item[\hskip\labelsep \theorem@headerfont ##1\ ##2.]}%
\def\@opargbegintheorem##1##2##3{%
   \item[\hskip\labelsep \theorem@headerfont ##1\ ##2\ (##3).]}}
\endgroup
\makeatother

\theoremstyle{upshape}

\newtheorem{defn}[theorem]{Definition}
\newtheorem{remark}[theorem]{Remark}

\makeatletter
\renewcommand{\subsection}{\@startsection{subsection}{2}{0pt}{-3ex
plus -1ex minus -0.2ex}{-2mm plus -0pt minus
-2pt}{\normalfont\bfseries}} 
\renewcommand{\subsubsection}{\@startsection{subsubsection}{3}{0pt}{-3ex
plus -1ex minus -0.2ex}{-2mm plus -0pt minus
-2pt}{\normalfont\bfseries}} 
\makeatother

\makeatletter
\@addtoreset{equation}{section}
\makeatother

\newcommand{\cntrct}                
{\hspace{2pt}\raisebox{1pt}{\text{$\lrcorner$}}\hspace{2pt}}

\newcommand{\proof}[1][Proof.]{\smallskip\noindent{\em #1}}
\def\endproof{\hfill\ensuremath{\square}\par\medskip}

\def\eqref#1{\thetag{\ref{#1}}}

\let\latexref=\ref
\def\ref#1{{\normalfont{\latexref{#1}}}}

\newcommand{\wt}{\widetilde}

\newcommand{\idot}{{\:\raisebox{1pt}{\text{\circle*{1.5}}}}}
\newcommand{\hdot}{{\:\raisebox{3pt}{\text{\circle*{1.5}}}}}
\renewcommand{\phi}{\varphi}

\def\dlim_#1{{\displaystyle\lim_{#1}}}

\newcommand{\Hom}{\operatorname{Hom}}
\newcommand{\Ext}{\operatorname{Ext}}

\newcommand{\Tor}{\operatorname{Tor}}

\newcommand{\Fun}{\operatorname{Fun}}

\newcommand{\id}{\operatorname{\sf id}}

\newcommand{\D}{{\cal D}}

\newcommand{\Sets}{\operatorname{Sets}}

\newcommand{\amod}{{\text{\rm -mod}}}

\newcommand{\ppt}{{\sf pt}}

\newcommand{\Ab}{\operatorname{Ab}}

\newcommand{\copr}{{\textstyle\coprod}}

\newcommand{\Z}{{\mathbb Z}}

\newcommand{\Top}{\operatorname{Top}}
\newcommand{\Real}{\operatorname{\sf Real}}

\newcommand{\vT}{t}

\newcommand{\U}{{\sf U}}
\newcommand{\T}{{\sf T}}

\newcommand{\Ll}{{\sf Q}}
\newcommand{\Rr}{{\sf R}}

\title{Homology of infinite loop spaces}

\author{D. Kaledin\thanks{Partially supported by AG Laboratory
    SU-HSE, RF government grant, ag. 11.G34.31.0023, the RFBR grant
    09-01-00242 and the Science Foundation of the SU-HSE award
    No. 10-09-0015}}

\begin{document}

\maketitle

\tableofcontents

\section*{Introduction.}

A {\em spectrum} $X_\idot$ is a sequence of pointed topological
spaces $X_0,X_1,\dots$ and homotopy equivalences $X_n \cong \Omega
X_{n+1}$, $n \geq 0$ (we tacitly assume that all the topological
spaces in consideration are nice enough, e.g. having homotopy type
of a CW complex). A spectrum $X_\idot$ is {\em connected} if all its
components $X_n$, $n \geq 0$ are connected. {\em Homology}
$H_\idot(X_\idot,\Z)$ of a spectrum $X_\idot$ with coefficients in a
ring $R$ is given by
$$
H_\idot(X_\idot,R) = \lim_{\overset{n}{\to}}\wt{H}_{\idot+n}(X_n,R),
$$
where $\wt{H}_\idot(-,R)$ denotes reduced homology of a pointed
topological space, and the limit is taken with respect to maps
$\Sigma X_n \to X_{n+1}$ adjoint to the structure maps $X_n \to
\Omega X_{n+1}$. For any $n,i \geq 0$, we then have a natural map
$$
\wt{H}_{i+n}(X_n,R) \to H_i(X_\idot,R).
$$
If the spectrum $X_\idot$ is connected, this map is an isomorphism
for $i < n$.

The forgetful functor sending a spectrum $X_\idot$ to its component
$X_0$ is conservative on the category of connected spectra, so that
up to a homotopy equivalence, a connected spectrum $X_\idot$ can be
reconstructed from a pointed topological space $X_0$ equipped with
an additional structure. This structure is usually called an
``infinite loop space structure'', and it can be described in
several ultimately equivalent ways, mostly discovered in the early
1970es and sometimes called ``machines'' (see \cite{adams} for an
all-time great overview of the subject). One of these machines is
that of G. Segal \cite{segal}, where a connected spectrum is
constructed from a so-called {\em special $\Gamma$-space}. This
turned out to be very useful, since e.g. in algebraic $K$-theory,
the relevant $\Gamma$-space often can be obtained almost for free.

\bigskip

The goal of this short note is to give a simple expression for the
homology of a connected spectrum $X_\idot$ in terms of the
associated special $\Gamma$-space. We state right away that the
expression is not new, and it is due to T. Pirashvili -- namely, it
can be deduced rather directly from \cite[Proposition 2.2]{pira3},
and for Eilenberg-Maclane spectra, the results goes back at least to
\cite{pira} (see the end of Section 3 for more details). All the
basic ideas behind the proof are also definitely due to
Pirashvili. However, the result itself is never stated explicitly in
the general corpus of Pirashvili's work, and while well-known to
experts, is not universally known. So, a short and self-contained
independent proof might be useful. This is what the present paper
aims to provide.

The paper consists of three parts: in Section 1, we recall the
details of the Segal machine in a convenient form, in Section 2, we
build a homological counterpart of the theory, and finally in
Section 3, we state and prove our results, and sketch an alternative
approach using \cite[Proposition 2.2]{pira3}.

\subsection*{Acknowledgements.} The paper owes its existence to
J. Peter May who explained to me that the result is not known to
everybody. I am also grateful to S. Prontsev for useful discussions,
and I am extremely grateful to T. Pirashvili who explained to me the
exact status of the result and kindly provided all the references. I
am grateful to the referee for useful suggestions.

\subsection*{A note on notation.} For the convenience of the reader,
here is a brief comparison between our notation and that of
Pirashvili. In \cite{pira3} and elsewhere, our $\Gamma_+$ is
$\Gamma$. Our functor $T$ is denoted $L$ in \cite{pira1}. Our $t$ is
$t$.

\section{Recollection on the Segal machine.}

We start by briefly recalling Segal's approach to infinite loop
spaces and rephrasing it in a language that suits our goal.

Denote by $\Gamma_+$ the category of finite pointed sets. For any
integer $n \geq 0$, denote by $[n]_+ \in \Gamma_+$ the set with $n$
unmarked elements (plus one distinguished element $o \in
[n]_+$). Alternatively, $\Gamma_+$ is equivalent to the category
$\Gamma'_+$ of finite sets and {\em partially defined maps} between
them -- that is, a map from $S_1$ to $S_2$ is given by a diagram
\begin{equation}\label{dom}
\begin{CD}
S_1 @<{\iota}<< S @>>> S_2
\end{CD}
\end{equation}
with injective $\iota$. The equivalence $\gamma:\Gamma_+' \to
\Gamma_+$ adds a distinguished element $o$ to a set $S \in
\Gamma'_+$, and for any $f:S_1 \to S_2$ represented by a diagram
\eqref{dom}, $\gamma(f)$ sends $S_1 \setminus \iota(S)$ to this
added distinguished element in $S_2$.

For any injective map $\iota:S_1 \to S_2$ of finite sets, let
$\iota^\#:S_2 \to S_1$ be the map in $\Gamma_+'$ represented by the
diagram
$$
\begin{CD}
S_2 @<{\iota}<< S_1 @>{\id}>> S_1.
\end{CD}
$$

\begin{defn}\label{spe.defn}
\begin{enumerate}
\item A {\em $\Gamma$-space} is a functor $X:\Gamma_+ \to \Top_+$
  from $\Gamma_+$ to the category of compactly generated pointed
  topological spaces. Say that a $\Gamma$-space $X$ is {\em
    normalized} if $X([0]_+)$ is the one-point set $\ppt$.
\item A $\Gamma$-space is {\em special} if it is normalized, and for
  any $S_1,S_2 \in \Gamma_+'$ with the natural embeddings
  $\iota_1:S_1 \to S_1 \copr S_2$, $\iota_2:S_2 \to S_1 \copr S_2$,
  the map
$$
\begin{CD}
X(\gamma(S_1 \copr S_2)) @>{X(\gamma(\iota_1^\#)) \times
  X(\gamma(\iota_2^\#))}>> X(\gamma(S_1)) \times X(\gamma(S_2))
\end{CD}
$$
is a homotopy equivalence.
\end{enumerate}
\end{defn}

\begin{remark}
Sometimes it is convenient to relax the normalization condition on
special $\Gamma$-spaces by only requiring that $X([0]_+)$ is
contractible. However, the stronger condition is harmless: replacing
$X([n]_+)$, $n \geq 1$, with the homotopy fiber of the map $X([n])_+
\to X([0]_+)$ corresponding to the unique map $[n]_+ \to [0]_+$, one
can always achieve $X([0]_+)=\ppt$.
\end{remark}

The category of $\Gamma$-spaces is denoted $\Gamma_+\Top_+$. For any
two $\Gamma$-spaces $X_1,X_2 \in \Gamma_+\Top_+$, we define $X_1
\vee X_2$, $X_1 \times X_2$ and $X_1 \wedge X_2$ pointwise.

\bigskip

Let $\wt{\Gamma_+\Top_+} \subset \Gamma_+\Top_+$ be the full
subcategory spanned by normalized $\Gamma$-spaces. Then the
forgetfull functor
$$
\U:\wt{\Gamma_+\Top_+} \to \Top_+, \qquad U(X) = X([1]_+)
$$
has a left-adjoint
$$
\T:\Top_+ \to \wt{\Gamma_+\Top_+}.
$$
Explicitly, for a pointed topological space $X$, the $\Gamma$-space
$\T(X)$ is given by
\begin{equation}\label{T.exp}
\T(X)([n]_+) = \bigvee_{s \in [n]_+ \setminus \{o\}}X = X \wedge
  [n]_+.
\end{equation}
The adjunction map $\id \to \U \circ \T$ is an isomorphism, so that
$\T$ is a full embedding, and the adjunction map $\tau:\T \circ \U
\to \id$ can be described as follows: for any $X \in
\Gamma_+\Top_+$, $[n]_+ \in \Gamma_+$, we have
$$
\tau = \bigvee_{s \in [n]_+ \setminus \{o\}}X(i_s):\T(\U(X))([n]_+) =
\bigvee_{s \in [n]_+ \setminus \{o\}}X([1]_+) \longrightarrow
X([n]_+),
$$
where $i_s:[1]_+ \to [n]_+$ is the embedding onto the subset
$\{s,o\} \subset [n]_+$.

\bigskip

Let $\Delta$ be, as usual, the category of finite non-empty totally
ordered sets, with $[n] \in \Delta$ denoting the set of integers
from $0$ to $n$, and let $S:\Delta^{opp} \to \Sets$ be the standard
simplicial circle -- that is, the simplicial set obtained by gluing
together the two ends of the standard $1$-simplex. The glued ends
give a natural distiguished element in $S([n])$, $[n] \in \Delta$,
so that $S$ is actually a pointed simplicial set. Moreover, $S([n])
\cong [n]_+$ is a finite set for any $[n] \in \Delta^{opp}$, so that
$S$ can be interpreted as a functor $\sigma:\Delta^{opp} \to
\Gamma_+$.

Recall that for any simplicial topological space $X:\Delta^{opp} \to
\Top_+$, we have its {\em geometric realization} $\Real(X) \in
\Top_+$, and this construction is functorial in $X$ and compatible
with products and colimits. For any simplicial abelian group $M$,
denote by $N_\idot(M)$ the corresponding standard complex with
terms $N_n(M) = M([n])$ and differential given by the alternating
sum of the face maps. Then for any ring $R$, the reduced singular
chain complex $\wt{C}_\idot(\Real(X),R)$ is naturally
quasiisomorphic to the total complex of a bicomplex
\begin{equation}\label{simpl.eq}
N_\idot(\wt{C}_\idot(X,R)).
\end{equation}
The geometric realization $\Real(S)$ of the simplicial circle is
homeomorphic to the $1$-sphere $S^1$. For any simplicial pointed
topological space $X:\Delta^{opp} \to \Top_+$, the realization
$$
\Real(X \wedge S)
$$
of the pointwise smash product $X \wedge S$ is homeomorphic to the
suspension $\Sigma X = S^1 \wedge X$.

\begin{defn}
{\em Geometric realization} $\Real(X)$ of a $\Gamma$-space $X$ is
given by
$$
\Real(X) = \Real(\sigma^*X).
$$
\end{defn}

In particular, for any $X \in \Top_+$, we have $\sigma^*\T(X) \cong
S \wedge X$, so that $\Real(\T(X)) \cong \Sigma X$.

\bigskip

Now consider the product $\Gamma_+^2 = \Gamma_+ \times
\Gamma_+$. Let $\pi_1,\pi_2:\Gamma_+^2 \to \Gamma_+$ be the
projections onto the first resp. second factor, and let
$\beta:\Gamma_+^2 \to \Gamma_+$ be the smash-product functor,
$\beta([n_1]_+ \times [n_2]_+) = [n_1]_+ \wedge [n_2]_+ \cong
[n_1n_2]_+$. Denote by $\Gamma_+^2\Top_+$ the category of functors
from $\Gamma_+^2$ to $\Top_+$, and let
$$
\U_1,\Real_1:\Gamma_+^2\Top_+ \to \Gamma_+\Top_+, \qquad
\T_1:\Gamma_+\Top_+ \to \Gamma_+^2\Top_+
$$
be the functors obtained by applying $\U$ resp. $\Real$ resp. $\T$
fiberwise over fibers of the projection $\pi_1:\Gamma_+^2 \to
\Gamma_+$ (in partricular, $\U_1 \cong i_1^*$, where $i_1:\Gamma_+ \to
\Gamma_+^2$ is the embedding onto $\Gamma_+ \times [1]_+ \subset
\Gamma_+^2$). For any normalized $\Gamma$-space $X:\Gamma_+ \to
\Top_+$, let
$$
BX = \Real_1(\beta^*X),
$$
and let
$$
\Sigma(X) = \Real_1(\T_1(X)).
$$
Note that for any $[n]_+$, we have 
$$
\Sigma(X)([n]_+) = \Real_1(\T_1(X))([n]_+) = \Real(\T(X([n]_+)))
\cong \Sigma X([n]_+),
$$
and in particular, $\U(\Sigma(X)) \cong \Sigma\U(X)$. Moreover,
since $\beta \circ i_1 \cong \id$, we have $\U_1(\beta^*X) \cong X$,
so that we obtain a natural adjunction map
\begin{equation}\label{tau.eq}
\tau:\T_1(X) \cong \T_1(\U_1(\beta^*X)) \to \beta^*X
\end{equation}
and its realization
\begin{equation}\label{rho.eq}
\rho_X = \Real_1(\tau):\Sigma(X) \to BX.
\end{equation}
Segal, then, proved the following.

\begin{prop}
Assume given a special $\Gamma$-space $X$. Then
\begin{enumerate}
\item the $\Gamma$-space $BX$ is also special, and
\item the natural map
$$
\U(X) \to \Omega\U(BX)
$$
adjoint to the map
\begin{equation}\label{trans.g}
\U(\rho_X):\Sigma \U(X) \cong \U(\Sigma(X)) \to \U(BX)
\end{equation}
is a homotopy equivalence.\endproof
\end{enumerate}
\end{prop}

By \thetag{i}, the functor $B$ can be iterated, so that every
special $\Gamma$-space $X$ gives rise to a sequence of special
$\Gamma$-spaces $B^nX$; by \thetag{ii}, the sequence $\U(B^nX)$ with
the maps $\U(\rho_{B^nX})$ then naturally forms a spectrum. We will
denote this spectrum by $EX_\idot$.

\section{Homology of $\Gamma$-spaces.}

Fix once and for all a commutative ring $R$, and consider the
category $\Fun(\Gamma_+,R)$ of functors from $\Gamma_+$ to the
category $R\amod$ of $R$-modules. This is an abelian category with
enough injectives and projectives. We equip it with pointwise tensor
product, and we denote by $\D(\Gamma_+,R)$ its derived category. An
obvious set of projective generators is given by representable
functors $T_n$,
$$
T_n([m]_+) = R(\Gamma_+([n]_+,[m]_+),
$$
since by Yoneda, we have $\Hom(T_n,E) \cong E([n]_+)$ for any $E \in
\Fun(\Gamma_+,R)$.  Let $T \in \Fun(\Gamma_+,R)$ be the functor
given by
$$
T([n]_+) = \overline{R[[n]_+]} = \bigoplus_{s \in [n]_+ \setminus
  \{o\}}R,
$$
that is, the reduced span functor. We have an obvious direct sum
decomposition $T_1 \cong T_0 \oplus T$.

Consider the functor $\T:R\amod \to \Fun(\Gamma_+,R)$ given by
$$
\T(M) = T \otimes_R M
$$
for any $R$-module $M$. This is consistent with previous notation,
in the sense that for any $X \in \Top_+$ with reduced singular chain
complex $\wt{C}_\idot(X,R)$, \eqref{T.exp} immediately gives a
canonical isomorphism
\begin{equation}\label{T.T}
\T(\wt{C}_\idot(X,R)) \cong \wt{C}_\idot(\T(X),R).
\end{equation}
The functor $\T:R\amod \to \Fun(\Gamma_+,R)$ is exact, and it has a
right and a left-adjoint $\Rr,\Ll:\Fun(\Gamma_+,R) \to R\amod$.

\begin{lemma}\label{R.T}
For any $E \in \Fun(\Gamma_+,R)$, we have a canonical decomposition
$E([1]_+) \cong R(E) \oplus E([0]_+)$. The functor $\Rr$ is exact,
the functor $\T$ is fully faithful, and its extension $\T:\D(R\amod)
\to \D(\Gamma_+,R)$ is also fully faithful.
\end{lemma}

\proof{} The decomposition is induced by the decomposition $T_1
\cong T \oplus T_0$. Exactness of $\Rr$ follows; to see that the
embedding $\T$ is fully faithful, note that $\Rr \circ \T \cong
\id$ both on the abelian and on the derived category level.
\endproof

\begin{defn}
The {\em homology} $H^\Gamma_\idot(E)$ of a functor $E \in
\Fun(\Gamma_+,R)$ is given by
$$
H_\idot^\Gamma(E) = L^\hdot\Ll(E),
$$
the derived functors of the functor $\Ll$ left-adjoint to the full
embedding $\T:R\amod \to \Fun(\Gamma_+,R)$.
\end{defn}

Explicitly, homology can be expressed as
$$
H^\Gamma_\idot(E) = \Tor_\idot^{\Gamma_+}(\vT,E),
$$
where $\vT:\Gamma^{opp}_+ \to R\amod$ is given by $\vT([n]_+) =
\Hom_R(T([n]_+),R)$, and $\Tor_\idot$ is taken over the small
category $\Gamma_+$ in the usual way, see e.g. \cite[Secton
  1.1]{ka.car}. To compute it, it suffices to find a projective
resolution of the functor $\vT$. One very elegant way to do it was
discovered by Pirashvili and Jibladze, and it leads to the so-called
cube construction of MacLane (see \cite{P}, or a slightly less
computational exposition in \cite[Section 3.3]{ka.car}). Whatever
resolution $Q_\idot$ one fixes, one immediately obtains a canonical
way to generalize homology to complexes: for any complex $E_\idot$
in $\Fun(\Gamma_+,R)$, we obtain a complex
\begin{equation}\label{Q.eq}
Q_\idot(E_\idot) = Q_\idot \otimes^{\Gamma_+} E_\idot
\end{equation}
of $R$-modules whose homology we denote by
$H^\Gamma_\idot(E_\idot)$.  If the complex $E_\idot = E$ is
concentrated in degree $0$, we have $H^\Gamma_\idot(E_\idot) \cong
H^\Gamma_\idot(E)$

The following Lemma is the crucial result of the theory.

\begin{lemma}\label{main.le}
Assume given $E_1,E_2 \in \Fun(\Gamma_+,R)$ such that $E_1([0]_+) =
E_2([0]_+) = 0$. Then $H^\Gamma_\idot(E_1 \otimes_R E_2)=0$.
\end{lemma}

\proof{} This is \cite[Lemma 2]{pira1}; I give a proof for the
convenience of the reader.

Consider the product $\Gamma_+ \times \Gamma_+$, let
$\pi_1,\pi_2:\Gamma_+ \times \Gamma_+ \to \Gamma_+$ be the
projections onto the first resp. second factor, and let
$\iota_1,\iota_2:\Gamma_+ \to \Gamma_+ \times \Gamma_+$ be the
embeddings sending $[n]_+$ to $[n]_+ \times [0]_+$ resp. $[0]_+
\times [n]_+$. Moreover, let $m:\Gamma_+ \times \Gamma_+ \to
\Gamma_+$ be the coproduct functor, and let $\delta:\Gamma_+ \to
\Gamma_+ \times \Gamma_+$ be the diagonal embedding. Then $\pi_i$ is
right-adjoint to $\iota_i$, $i=1,2$, and $m$ is right-adjoint to
$\delta$. We obviously have
$$
m^*T \cong \pi_1^*T \oplus \pi_2^*T,
$$
and this decomposition induces an isomorphism
$$
m^* \circ \T \cong (\pi_1^* \circ T) \oplus (\pi_2^* \circ T).
$$
By adjunction, we obtain a functorial isomorphism
$$
L^\hdot\Ll(\delta^*E) \cong L^\hdot\Ll(\iota_1^*E) \oplus
L^\hdot\Ll(\iota_2^*E)
$$
for any functor $E:\Gamma_+ \times \Gamma_+ \to R\amod$. Take $E=E_1
\boxtimes_R E_2$, and note that by assumption, $\iota_1^*E =
\iota_2^*E = 0$, while $\delta^*E \cong E_1 \otimes_R E_2$.
\endproof

\begin{defn}\label{ho.ga.de}
The homology $H_\idot^\Gamma(X,R)$ of a $\Gamma$-space $X$ is given
by
$$
H_\idot^\Gamma(X,R) = H_\idot^\Gamma(\wt{C}_\idot(X,R)),
$$
where $\wt{C}_\idot(X,R)$ is a complex in $\Fun(\Gamma_+,R)$
obtained by taking pointwise the reduced singular chain complex
$\wt{C}_\idot(-,R)$.
\end{defn}

\begin{lemma}\label{T.H}
For any $X \in \Top_+$, we have a canonical isomorphism
$$
\wt{H}_\idot(X,R) \cong H^\Gamma_\idot(\T(X),R).
$$
\end{lemma}

\proof{} By virtue of the quasiisomorphism \eqref{T.T}, this
immediately follows from the last claim of Lemma~\ref{R.T}: we have
$$
\Tor^{\Gamma_+}_\idot(t,T) \cong
\Hom_R(\Ext^\hdot_{\Gamma_+}(T,T),R) \cong
\Hom_R(\Ext^\hdot_R(R,R)),
$$
and the right-hand side is $R$ in degree $0$ and $0$ in higher
degrees, so that for any complex $C_\idot$ of $R$-modules, the
groups 
$$
H^\Gamma_\idot(\T(C_\idot)) \cong \Tor_\idot^{\Gamma_+}(t,C_\idot
\otimes T) 
$$
coincide with the homology groups of the complex $C_\idot$ itself.
\endproof

Lemma~\ref{main.le} has the following implication for the homology
of $\Gamma$-spaces. For any $n \geq 0$, let $\mu_n:\Gamma_+ \to
\Gamma_+$ be the functor given by $\mu_n([m]_+) = [n]_+ \wedge
[m]_+$.  Assume given two pointed finite sets $[n]_+,[n']_+ \in
\Gamma_+$, identify $[n]_+ \vee [n']_+ \cong [n+n']_+$, and consider
the natural maps $[n]_+ \to [n]_+ \vee [n']_+ \cong [n+n']_+$,
$[n']_+ \to [n]_+ \vee [n']_+ \cong [n+n']_+$. These maps then
induce maps
$$
\iota:\mu_n \to \mu_{n+n'}, \qquad \iota':\mu_{n'} \to \mu_{n+n'},
$$
and for any $\Gamma$-space $X$, we obtain natural maps
$$
\iota:\mu_n^*X \to \mu_{n+n'}^*X, \qquad \iota':\mu_{n'}^*X \to
\mu_{n+n'}^*X
$$
and
$$
\iota^{\#}:\mu_{n+n'}^*X \to \mu_n^*X, \qquad
\iota^{'\#}:\mu_{n+n'}^*X \to \mu_{n'}^*X.
$$

\begin{corr}\label{Coro}
Assume that the $\Gamma$-space $X$ is special, Then the natural map
$$
\iota \vee \iota':\mu_n^*X \vee \mu_{n'}^*X \to \mu_{n+n'}^*X
$$
of $\Gamma$-spaces induces an isomorphism of homology
$H_\idot^\Gamma(-,R)$.
\end{corr}

\proof{} For any two pointed topological spaces $X_1$, $X_2$, we
have a cofiber sequence
$$
\begin{CD}
X_1 \vee X_2 @>>> X_1 \times X_2 @>>> X_1 \wedge X_2.
\end{CD}
$$
Since $X$ is special, the natural map
$$
\begin{CD}
\mu_{n+n'}^*X @>{\iota^{\#} \times \iota^{'\#}}>> \mu_n^*X \times
\mu_{n'}^*X
\end{CD}
$$
is a pointwise homotopy equivalence. Therefore the sequence
$$
\begin{CD}
\mu_n^*X \vee \mu_{n'}^*X @>{\iota \vee \iota'}>> \mu_{n+n'}^*X
@>{\iota^{\#} \wedge \iota^{'\#}}>> \mu_n^*X \wedge \mu_{n'}^*X
\end{CD}
$$
is a pointwise cofiber sequence, and it suffices to prove that
$$
H_\idot^\Gamma(\mu_n^*X \wedge \mu_{n'}^*X,R)=0.
$$
This immediately follows from Lemma~\ref{main.le} and the K\"unneth
formula.
\endproof

\section{Stabilization.}

We can now formulate and prove the main result of the paper. For any
special $\Gamma$-space $X$, let
$$
\tau_X:\T(\U(X)) \to X
$$
be the adjunction map, and let $\rho_X:\Sigma(X) \to BX$ be as in
\eqref{rho.eq}.

\begin{lemma}
For any special $\Gamma$-space $X$, the diagram
$$
\begin{CD}
\Sigma(\T(\U(X))) \cong \T(\U(\Sigma(X))) @>{\Sigma(\tau_X)}>>
\Sigma(X)\\
@V{\T(\U(\rho_X))}VV @VV{\rho_X}V\\
\T(\U(BX)) @>{\tau_{BX}}>> BX
\end{CD}
$$
is commutative.
\end{lemma}

\proof{} By \eqref{T.exp}, we have
$$
\begin{aligned}
\T_1(\T(Y))([n_1]_+ \times [n_2]_+) &\cong [n_2]_+ \wedge
\T(Y)([n_1]_+) \cong [n_2]_+ \wedge [n_1]_+ \wedge Y\\
&\cong [n_1n_2]_+ \wedge Y \cong \beta^*\T(Y)([n_1]_+ \times
  [n_2]_+)
\end{aligned}
$$
for any $Y \in \Top_+$, $[n_1]_+,[n_2]_+ \in \Gamma_+$, so that
$\T_1(\T(Y)) \cong \beta^*\T(Y)$. Taking $Y = \U(X)$, we obtain a
natural commutative diagram
$$
\begin{CD}
\T_1(\T(\U(X)) @>>> \T_1(X)\\
@| @VVV\\
\beta^*\T(\U(X)) @>>> \beta^*X.
\end{CD}
$$
Applying $\Real_1$, we get the claim.
\endproof

Taking homology $H^\Gamma_\idot(-,R)$ and using Lemma~\ref{T.H}, we
obtain a commutative diagram
$$
\begin{CD}
\wt{H}_\idot(\U(X),R) @>>> H^\Gamma_\idot(X,R)\\
@VVV @VVV\\
\wt{H}_{\idot+1}(\U(BX),R) @>>> H^\Gamma_{\idot+1}(BX,R),
\end{CD}
$$
and passing to the limit, we get a natural map
\begin{equation}\label{main.eq}
H_\idot(EX_\idot,R) \to
\lim_{\overset{n}{\to}}H^\Gamma_{\idot+n}(\U(B^nX),R).
\end{equation}
Here is, then, our main result.

\begin{theorem}\label{main}
Assume given a special $\Gamma$-space $X$, and let $EX_\idot$ be the
corresponding spectrum. Then the natural map \eqref{main.eq} factors
through an isomorphism
$$
H_\idot(EX_\idot,R) \cong H_\idot^\Gamma(X,R),
$$
where the right-hand side is as in Definition~\ref{ho.ga.de}.
\end{theorem}

The proof is a combination of the following two results.

\begin{lemma}\label{1.le}
Assume given a special $\Gamma$-space $X$. Then the map $\rho_X$ of
\eqref{rho.eq} induces an isomorphism
$$
H^\Gamma_\idot(\Sigma(X),R) \cong H^\Gamma_\idot(BX,R).
$$
\end{lemma}

\proof{} Combining \eqref{simpl.eq} and \eqref{Q.eq}, we see that
for any $X \in \Gamma^2_+\Top_+$, the homology
$H^\Gamma_\idot(\Real_1(X),R)$ can be computed by the total complex
of the triple complex $Q_\idot(N_\idot(C_\idot(X),R))$. This gives
rise to a convergent spectral sequence
$$
H^\Gamma_i(i_n^*X,R) \Rightarrow H^\Gamma_{i+n}(\Real_1(X),R),
$$
where $i_n:\Gamma_+ \to \Gamma^2_+$ is the embedding onto $\Gamma_+
\times [n]_+$. We conclude that to prove the lemma, it suffices to
prove that for every $n \geq 0$, the map
$$
H^\Gamma_\idot(i_n^*\T_1(X),R) \to
H^\Gamma_\idot(i_n^*\beta^*X,R)
$$
induced by the map $\tau$ of \eqref{tau.eq} is an isomorphism. By
\eqref{T.exp}, we have
$$
i_n^*\T_1(X) \cong [n]_+ \wedge X \cong \bigvee_{s \in [n]_+
  \setminus \{o\}}X
$$
and by definition,
$$
i_n^*\beta^*X \cong \mu_n^*X.
$$
so that the statement immediately follows by induction on $n$ from
Corollary~\ref{Coro}.
\endproof

\begin{lemma}\label{2.le}
Assume given a special $\Gamma$-space $X$, and assume that $\U(X)$
is $n$-connected for some $n \geq 1$. Then the natural map
$$
\wt{H}^\Gamma_i(\U(X),R) \cong H^\Gamma_i(\T(\U(X)),R) \to
H^\Gamma_i(X,R)
$$
induces by the map $\tau_X$ is an isomorphism for $i < 2n$.
\end{lemma}

\proof{} Since $\U(X)$ is $n$-connected, $H_0(\U(X),R) \cong R$ and
$H_i(\U(X),R) = 0$ when $0 < i < n$. Then by
Definition~\ref{spe.defn}~\thetag{ii}, $X([n]_+)$ is
homotopy-equivalent to $X([1]_+)^n = \U(X)^n$ for any $n \geq 1$,
and then the K\"unneth formula immediately implies that the map
$$
\wt{C}_\idot(\T(\U(X)),R) \to \wt{C}_\idot(X,R)
$$
induced by $\tau_X$ is a quasiisomorphism in homological degrees
less than $2n$. Applying $L^\hdot\Ll$, we get the claim.
\endproof

\proof[Proof of Theorem~\ref{main}.] By Lemma~\ref{1.le}, the
natural map
$$
H^\Gamma_\idot(X,R) \to
\lim_{\overset{n}{\to}}H^\Gamma_{\idot+n}(\U(B^nX),R)
$$
is an isomorphism, and by Lemma~\ref{2.le}, \eqref{main.eq} is also
an isomorphism -- in fact, in each homological degree, it becomes an
isomorphism at some finite step in the inductive sequence.
\endproof

To finish the paper, let us explain how Theorem~\ref{main} can be
deduced from the work of T. Pirashvili mentioned in the
introduction. Note that any abelian group can be treated as a
pointed set, by taking $0$ as the distinguished point and forgetting
the rst of the group structure. Thus a functor $F:\Gamma_+ \to
\Delta^{opp}\Ab$ from $\Gamma_+$ to the category of simplicial
abelian groups can be treated as a pointed simplicial $\Gamma$-set.
Then even if such $F$ is not special in the sense of
Definition~\ref{spe.defn}~\thetag{ii}, the map \eqref{trans.g} is
still well-defined, so that the sequence $B^nF$, $n \geq 0$ forms a
pre-spectrum. One denotes by $\pi_\idot^{st}(F)$ the homotopy groups
of the corresponding spectrum. Then \cite[Proposition 2.2]{pira3}
claims that there exists a natural isomorphism
$$
\pi_\idot^{st}(F) \cong \Tor_\idot^{\Gamma_+}(t,F)
$$
(to be precise, \cite[Proposition 2.2]{pira3} is stated only for
functors to constant simplicial groups, but generalization to
arbitrary ones is immediate). Pirashvili's proof of this fact also
uses Lemma~\ref{main.le}, but it is in fact simpler since working
with an arbitrary $F$ gives more lattitude. Then to deduce
Theorem~\ref{main}, one has to take $F = C_\idot(X,R)$, and show
that
$$
H_\idot(EX_\idot,R) \cong \pi_\idot^{st}(F).
$$
This is also rather straighforward. So in a nutshell, Pirashvili's
proof is ultimately simpler but relies on some context, while our
proof is longer but elementary and self-contained.

\medskip

\noindent
{\sc
Steklov Math Institute\\
Moscow, USSR}

\medskip

\noindent
{\em E-mail address\/}: {\tt kaledin@mi.ras.ru}

\end{document}